\documentclass[12pt]{amsart}
\usepackage{amsfonts}
\newtheorem{corollary}{Corollary}
\newtheorem{theorem}{Theorem}
\newtheorem{lemma}{Lemma}
\newtheorem{proposition}{Proposition}

\newtheorem{definition}{Definition}
\newtheorem{claim}{Claim}

\newcommand{\bee}[1]{\begin{equation}\label{#1}}
\newcommand{\beq}[1]{\begin{eqnarray}\label{#1}}
\newcommand{\ene}{\end{equation}}
\newcommand{\eqe}{\end{eqnarray}}
\newcommand{\ld}{\ldots}

\newcommand{\wg}{\widehat{G}}

\newcommand{\wG}{\widetilde{G}}
\newcommand{\wP}{\widetilde{P}}
\newcommand{\vp}{\varphi}

\newcommand{\ve}{\varepsilon}

\newcommand{\bG}{\overline{G}}

\newcommand{\bP}{\overline{P}}

\newcommand{\bM}{\overline{M}}
\newcommand{\bN}{\overline{N}}

\newcommand{\Sp}[1]{\mathrm{Span}\,\{#1\}}
\newcommand{\da}[2]{\mathrm{id}^p_{#2}\{#1\}}
\newcommand{\sa}[1]{\mathrm{alg}^p\{#1\}}

\newcommand{\rl}{restricted Lie algebra }
\newcommand{\rlc}{restricted Lie algebra, }
\newcommand{\rlp}{restricted Lie algebra. }
\newcommand{\rls}{restricted Lie algebras }
\newcommand{\rlsc}{restricted Lie algebras, }
\newcommand{\rlsp}{restricted Lie algebras. }
\newcommand{\fr}{free restricted Lie algebra }
\newcommand{\frs}{free restricted Lie algebras }
\newcommand{\frc}{free restricted Lie algebra, }
\newcommand{\frp}{free restricted Lie algebra. }
\newcommand{\ad}{\mathrm{ad}\,}
\newcommand{\nat}{\mathbb{N}}

\newcommand{\de}[2]{\delta_{#1}(#2)}
\newcommand{\gp}[1]{g^{[p^{#1}]}}
\newcommand{\hp}[1]{h^{[p^{#1}]}}
\newcommand{\tp}[1]{t^{[p^{#1}]}}
\begin{document}

\title{Large Restricted Lie Algebras}
\author{Yuri Bahturin}\address{Department of Mathematics and Statistics\\Memorial University of Newfoundland\\St. John's, NL, A1C 5S7, \textsc{Canada} and Department of Algebra\\Faculty of Mathematics and Mechanics\\119899 Moscow, \textsc{Russia}}\thanks{$^1$ Partially supported by NSERC grant \# 227060-04 and URP grant, Memorial University of Newfoundland}
\author{Alexander Olshanskii}\address{Department of Mathematics\\
1326 Stevenson Center\\
Vanderbilt University\\ and Department of Algebra\\Faculty of Mathematics and Mechanics\\119899 Moscow, \textsc{Russia}}\thanks{$^2$ Partially supported by NSF grants DMS-0245600 and DMS-0455881 and by RFBR grant  05-01-00895}
\begin{abstract} We establish some results about large \rls similar to those known in the Group Theory. As an application, we use this group-theoretic approach to produce some examples of restricted as well as ordinary Lie algebras which can serve as counterexamples for various Burnside - type questions.
\end{abstract}
\maketitle
\section{Introduction}\label{sI}

In this paper we consider \textit{restricted} Lie algebras over a \textit{perfect} field $F$ of characteristic $p>0$. For the basic information see \cite{J} and \cite{B}. All subalgebras and ideals considered will also be restricted, that is, closed under the $p$-operation $x\mapsto x^{[p]}$. 

\begin{definition}
A \rl $G$ is called \emph{large} if there is a subalgebra $H$ in $L$ of finite codimension such that $H$ admits a surjective homomorphism on a nonabelian \frp
\end{definition}
 
The main goal of this paper is to prove three theorems about restricted Lie algebras presented in terms of generators and defining relations. The first two deal with the large \rls while the third one is an application of the methods used in the first two theorems to the construction of finitely generated nil \rls of infinite dimension and their generalizations. A \rl $G$ is called \textit{nil} if for any $g\in G$ there is natural $n$ such that $\gp{n}=0$.

We introduce some more notation. Given a set $X$ of elements of a \rl $G$ we denote by $\sa{X}$ the restricted subagebra $H$ of $G$ generated by $X$ (if $X=\{ g\}$ then we simply write $\sa{g}$). Any element of $H$ has the form of $\sum_{i=1}^m [f_i(x_1,\ldots,x_n)]^{[p^{k_i}]}$ where each $f_i(x_1,\ldots,x_n)$ is an ordinary Lie polynomial in $x_1,\ldots,x_n\in X$. Another notation, $\da{X}{G}$, will be used to denote the restricted ideal $I$ of $G$ generated by $X$. Again, each element of $I$ will look like $\sum_{i=1}^m w_i^{[p^{k_i}]}$ where each $w_i$ is in the ordinary ideal of $G$ generated by $X$. While obtaining these remarks, it is important to remember the main identity of \rls: 
\begin{equation}\label{em}
[\gp{},h]=(\ad g)^p(h)=[\underbrace{g,\ldots,g}_{p},h].
\end{equation}

We recall that given $g\in G$ and $n$ a nonnegative integer, one defines $\gp{n}$ by induction if one sets $\gp{0}=g$ and $\gp{n}=(\gp{n-1})^{[p]}$ if $n\ge 1$. We call $g$ \textit{nilpotent} if $\gp{n}=0$ for some $n\in\nat$. In this case, by (\ref{em}), we also have that the linear transformation $\ad g :G\rightarrow G$ defined by $(\ad g)(x)=[g,x]$ for any $x\in G$ is also nilpotent. A nonzero element $g\in G$ is called \textit{algebraic} if $\dim(\sa{g})=n<\infty$. In this case $G=\Sp{g,\gp{},\ldots,\gp{n-1}}$, for some $n\in\nat$. If, additionally, $g$ is nilpotent then  $p^n$ is called the \textit{nil-index} of $g$. Thus $n$ is the least natural number such that $\gp{n}=0$. A \rl $G$ is called \textit{cyclic} if $G=\sa{g}$ for some $g\in G$.  

\medskip

Our first result is an analogue of  a group-theoretic theorem due to B. Baumslag and S. Pride \cite{BP}.

\begin{theorem}\label{tbp} Let $G$ be a restricted Lie algebra over a perfect field $F$ of characteristic $p >0$ given by a presentation with $n$ generators and $m$ relations, where  $m\leq n-2$. Then $G$ is large.
\end{theorem}

A technically useful form of Theorem \ref{tbp}, immediate from its proof, is as follows.

\begin{proposition}\label{r001}   Let $G$ be a restricted Lie algebra over a perfect field $F$ of characteristic $p >0$, given by a presentation with $n$ generators and $m$ relations, where  $m\leq n-2$. Then for any cyclic \rl $A$ of sufficiently large dimension, $G$ has a restricted ideal $M$, with $G/M\cong A$, such that $M$ maps homomorphically on a nonabelian \frp
\end{proposition}

The next theorem is an analogue of a group-theoretic theorem due to M. Lackenby \cite{L}, with a simplified proof due to A. Ol'shanskii - D. Osin \cite{OO}.

\begin{theorem}\label{tLOO} Let $G$ be a restricted Lie algebra over a \textit{perfect} field $F$ of characteristic $p >0$, $H$ an ideal of finite codimension in $G$ admitting a homomorphism on a nonabelian free restricted Lie algebra, $g_1,\ldots,g_k$ a set of elements of $H$. Let $I_n$ be a (restricted) ideal of $G$ generated $g_1^{[p^n]},\ldots,g_k^{[p^n]}$. Then $G/I_n$ is large for all but finitely many $n\in\mathbb{N}$.
\end{theorem}

An important particular case of this theorem, with some more information, reads like this.

\begin{proposition}\label{p881} Let $L$ be a \fr of rank at least two, $N$ an ideal of finite codimension in $L$, $g\in L\setminus N$ and $h\in\sa{g}\cap N$. Then $N/\da{h}{L}$ is a large \rl with a presentation in which the number of generators exceeds the number of defining relations at least by 2. As a consequence, also $L/\da{h}{L}$ is a large \rlp
\end{proposition}

The above results allow us to construct some examples in the spirit of the Unrestricted Burnside Problem for groups. In \cite[Chapter V, Exercise 17]{J} the author asked for the proof of the finite-dimensionality (probably under certain conditions) of finitely generated nil \rlsp Examples of infinite - dimensional finitely generated nil  \rls can be derived from E. Golod's original example of finitely generated Engel Lie algebras \cite{G}. In distinction with the situation in the Group Theory, where the example giving negative solutions to the Unrestricted Burnside problem are abundant, in the case of Lie algebras until now we had just one Golod's example and its derivatives.

Before we formulate our results, we recall some terminology.

An algebra is called \textit{residually finite-dimensional} (respectively, \textit{residually finite-dimensional nilpotent}) if for every nonzero element $g\in G$ there is a homomorphism $\vp$ of $G$ onto a finite-dimensional (respectively, finite-dimensional nilpotent) \rl such that $\vp(g)\neq0$. Equivalently, one can say that $G$ is residually finite-dimensional (nilpotent) if $G$ has a set of ideals $\{ I_{\alpha}\}$ with trivial intersection $\cap_{\alpha}I_{\alpha}$ and such that each quotient algebra $G/I_{\alpha}$ is finite-dimensional (nilpotent). The same definition applies to ordinary \rlsp A \textit{subfactor} of a \rl $G$ is a \rl $H/K$ where $K$ is an ideal of $H$ and $H$ is a subalgebra of $G$.

\begin{theorem}\label{tOO} 
Let $F$ be a \textit{perfect} at most countable field $F$ of characteristic $p >0$. Then for any finitely generated \rl $G$ with an ideal $P$ of finite codimension that can be mapped homomorphically onto a nonabelian \fr there exists an infinite-dimensional homomorphic image $\wG$ in which the image $\wP$ of $P$ is a nil \rlc and $\wG/\wP\cong G/P$. One can choose $\wG$ residually finite - dimensional and an inductive limit of large \rlsp
\end{theorem}

\begin{corollary}\label{cOO0} 
Let $F$ be a \textit{perfect} at most countable field of characteristic $p >0$. Then any \rl $G$, with a presentation where the number of generators exceeds the number of relations at least by two, there exists an infinite-dimensional homomorphic image $\wG$  which is a nil \rlp One can choose $\wG$ residually finite - dimensional nilpotent and an inductive limit of large \rlsp
\end{corollary}

\begin{corollary}\label{c001} 
Let $F$ be a \textit{perfect}  at most countable field of characteristic $p >0$. Then there exist infinite-dimensional finitely generated nil \rlsp  One can choose such algebras residually finite - dimensional, inductive limits of large \rlsc and with each finite dimensional subfactor nilpotent.
\end{corollary}

\begin{corollary}\label{c002} 
\rmfamily{\upshape{(E.S. Golod)}} Let $F$ be an  at most countable field of characteristic $p >0$. Then there exist infinite-dimensional finitely generated Engel Lie algebras. One can choose such algebras residually finite - dimensional nilpotent.
\end{corollary}

\bigskip

\section{Some Properties of Large Restricted Lie Algebras}

First, we want to mention a couple of obvious properties of large \rlsc following from the additivity of codimenesion.

\begin{proposition}\label{ppla1} The following are true.
\begin{enumerate}
\item[(i)] If a \rl $G$ has a homomorphic image which is large then also $G$ is large;
\item[(ii)] If a subalgebra of finite codimension in a \rl $G$ is large then also $G$ is large. 
\end{enumerate}
\end{proposition}

The next result requires  a little more sophistication.

\begin{proposition}\label{ppla2} The following are true.
\begin{enumerate}
\item[(i)] If a \rl can be mapped onto a nonabelian \fr then any subalgebra of finite codimension has the same property. A subalgebra of finite codimension in a large \rl is itself large.
\item[(ii)] If a subalgebra $H$ of finite codimension in a \rl $G$ can be mapped onto a nonabelian \fr then an ideal $K$ of finite codimension in $G$ also has this property. One can choose $K$ with $K\subset H$.
\end{enumerate}
\end{proposition}

\begin{proof}  To prove (i), we notice that if a \rl $G$ can be mapped onto a \fr $L$ by means of a surjective homomorphism $\ve$ and $H$ is a subalgebra of finite codimension in $G$ then $M=\ve(H)$ is a subalgebra of finite codimension in $M$. If $r$ is the number of generators in $L$ (could be an infinite cardinal) and $d=\dim L/M$ then an analogue of Schreier's formula for groups, due to G.P.Kukin, see for example \cite[2.7.5]{B}, says that the number of free generators for $M$ is given by $p^d(r-1)+1$. Obviously, if $r$ is greater than 1, this latter number is greater than 1, proving that $M$ is indeed a nonabelian \frp 

The second claim in (i) now follows since if a \rl $G$ is large, $G_1$ a restricted subalgebra of $G$ of finite codimension and a restricted subalgebra $H$ of finite codimension in $G$ can be mapped on a \fr then $G_1\cap H$ is a restricted subalgebra of finite codimension in $G_1$ which by what we have just proved can be mapped onto a \frp Thus $G_1$ is large.

Before we prove Claim (ii), we need a very general module-theoretic result.

\begin{lemma}\label{lrm}
Let $S$ be a unital subalgebra of an associative algebra $R$ with $1$ over a field $F$,  such that $R$ is generated as a left regular $S$-module by a finite subset $T$. Assume that $U$ is a unital left $R$-module, $V$ an $S$-submodule of $V$ such that $\dim_F U/V<\infty$. Then there is an $R$-submodule $W$ such that $W\subset V$ and still $\dim_F U/W<\infty$.
\end{lemma}
 
\begin{proof} For each $r\in R$ we consider a linear mapping $\Phi(r): V\rightarrow U/V$ given by $\Phi(r)(v)=rv+V$, for any $r\in R$ and $v\in V$. Obviously, $\Phi(s)(V)=sV=\{ V\}$, for any $s\in S$. The set $W=\cap_{r\in R}\mathrm{Ker}\,\Phi(r)$ is is easily seen to be an $R$-submodule of $U$ contained in $V$. Also, if $\Phi(t)(v)=V$ for all $t\in T$ then $v\in W$. This follows because for any $s\in S$, any $t\in T$ and $v$ as just above we have $\Phi(st)(v)=(st)v+V=s(tv)+V=V$. As a result, $W$ contains the intersections of the kernels of the finite set of linear mappings $ \Phi(t)$, $t\in T$, into a finite-dimensional space $U/V$. Each such kernel is of finite codimension by the Isomorphism Theorem, proving that $\dim V/W<\infty$.
\end{proof}

Now we can continue with the proof of Claim (ii). A subspace $V$ of a \rl $G$ is a restricted ideal of $G$ if and only if $V$ satisfies two conditions. First, $V$ must be a submodule under the natural $R$-module structure of $G$, $R$ the restricted enveloping algebra of $G$. Second, $V$ must be closed under the $p$-operation of $G$. If $G$ has a restricted subalgebra $H$ of finite codimension then by PBW-Theorem \cite[Chapter 5]{J} $R$ is a (free) finitely generated left (and right!) module over the associative subalgebra $S$ generated by $H$. Now Lemma \ref{lrm} with $U=G$ and $V=H$ applies and provides us with a subspace $W$ of finite codimension in $H$. Since $H$ is closed under the $p$-map, the $p$-closure $K$ of $W$ is a restricted ideal of $G$ contained between $W$ and $H$. By Part (i) of this Lemma it follows that $K$ can be mapped onto a \frc as required.
\end{proof}

In view of the last result one can define a large \rl as one with an \textit{ideal} of finite codimension which can be homomorphically mapped onto a nonabelian \frp

\section{Baumslag - Pride's Theorem for Restricted Lie Algebras}\label{sBPT}

Let $F$ be a perfect field of characteristic $p>0$ and $L=L(X)$ a free restricted Lie algebra over $F$ with a set of free generators $X=\{ x_1,\ldots,x_n\}$. Let also $W=\{ w_1,\ldots,w_m\}$ be a set of elements in $L$, $I$ a  restricted ideal in $L$ generated by $W$, and $G=L/I$. We then say that $G$ has a presentation $G=\langle x_1,\ldots,x_n\,|\,w_1,\ldots,w_m\rangle$ with $n$ generators and $m$ relations $w_1=0,\ld,w_m=0$. Sometimes the left hand sides of the relations, that is, the elements of $W$ are called the \textit{relators}. In \cite{BP} it was established that a \textit{group} which can be presented by $n$ generators and $m$ defining relations is large provided that $m\leq n-2$. We want to adapt this result to our situation.

Before we formulate our first result, we recall \cite{B} that in a free restricted Lie algebra $L(X)$ any element $w$ can be uniquely written as a linear combination of $p$-powers of generators (the power component) plus a linear combination of commutators of degree at least two in the generators and their $p$-powers (the commutator component).

The first result we would like to start with is the following.

\begin{proposition}\label{p001}  Let a restricted Lie algebra $G$ be presented in terms of generators and defining relations as above, with $m\leq n-1$. Then another presentation can be chosen so that one of the generators is not involved in the power components of the defining relations.
\end{proposition}

\begin{proof} Given a set of elements $a_1,\ld,a_d$ of a \rl $G$, the following transformation is called an \textit{elementary transformation}: 
$a_i\rightarrow \lambda a_i+f(a_1,\ld,a_{i-1},a_{i+1},\ld,a_d), a_k\rightarrow a_k$, for some $i\neq j$, and any $k\neq i,j$, $f$ a $p$-polynomial, $\lambda$ a nonzero element of $F$.

It is well known \cite[Chapter 2]{B} that any set $b_1,\ld,b_d$ that can be obtained by a finite sequence of these transformations from $a_1,\ld,a_d$ generates the same restricted subalgebra of $G$. In the case where $a_1,\ld,a_d$ is the set of free generators of a \fr $G$ the set $b_1,\ld,b_d$ is again the set of free generators of $G$. Each elementary transformation extends to an automorphism of $G$ and it is known that the group of automorphisms of a \fr is generates by such automorphisms.

We need an easy result about free abelian restricted Lie algebras.

\begin{lemma}\label{l001}Suppose we are given a free abelian restricted Lie algebra $A$ with free generators $y_1,\ldots,y_n$ and a restricted subalgebra $B$ generated by a set of elements $v_1,\ldots,v_m$. Then there is another free generating set $z_1,\ldots,z_n$ for $A$ and a set of elements\ \ $w_1,\ldots,w_k$ generating $B$ such that $k\leq \min(m,n)$,  and each $w_i$ is a $p$-polynomial in $z_i$, $i=1,\ldots,k$.
\end{lemma}

\begin{proof} It is well-known from \cite{J} that any abelian restricted Lie algebra $A$ is a left module over a skew polynomial algebra $\Lambda$ in one variable $t$ in the sense that $\alpha^p t=t\alpha$, for any $\alpha\in F$. The action is given by $t\ast a= a^{[p]}$. The left and right analogues of the Division Algorithm work in $\Lambda$ provided that the base field is perfect. If $A$ is a free abelian restricted Lie algebra as above then $A$ is a free left $\Lambda$-module with the free generating set $y_1,\ldots,y_n$. The $\Lambda$-submodules of $A$ are precisely the subalgebras of $A$.

For the proof of our lemma we need to write the matrix $(f_{ij}(t))$ of coefficients of the expression of $v_i$ in terms of $y_j$. As in the case of the Fundamental Theorem of Finitely Generated Modules Over Principal Ideal Domains, we apply elementary transformations to the rows and columns of this matrix. The elementary transformations of the rows of this matrix correspond to elementary transformations of $u_1,\ld,u_m$, which replace one generating set of $B$ by another, in particular, replacing $v_i$ by $v_i+g(t)\ast v_j$, for $i\neq j$, corresponds to adding to the $i^{th}$ row the $j^{th}$ one multiplied on the left by $g(t)$. Elementary transformations of the columns correspond to the replacement of one free generating set of $A$ by another. For example, if we modify the $j^{th}$ column by subtracting from it the $i^{th}$ one, $i\ne j$, then we obtain the matrix of coefficients of the generators of $B$ with respect to the new free generating set where $y_i$ is replaced by $y_i+g(t)\ast y_j$ while the remaining ones are left intact.

If we apply the natural versions of the left and the right division algorithms in $\Lambda$ to the above matrix it becomes obvious that using the elementary transformations we can reduce our matrix to the form where the only nonzero elements are the first $k$ diagonal elements $f_{11}(t)$, \ldots, $f_{kk}(t)$ where $k\leq \min(m,n)$.
\end{proof}

Now we can continue with the proof of  Proposition \ref{p001}. Let us assume that $L$ is freely generated by $x_1,\ld,x_n$ and $J$ the restricted ideal of $L$ generated by $u_1,\ld,u_m$ so that $G=L/J$. Recall that $m\leq n-1$. We can work modulo the commutator subalgebra $[L,L]$ of $L$. Suppose $A=L/[L,L]$ is the respective \frp Let $y_1,\ld,y_n$, $v_1,\ld,v_m$, and $B$ be the images of $x_1,\ld,x_n$, $u_1,\ld,u_m$, and $J$ under the natural homomorphism of $L$ onto $A$. Then we apply Lemma \ref{l001}. As a result, we obtain the free generators $z_1,\ld,z_n$ of $A$ and the generators $w_1,\ld,w_k$ of $B$, each $w_i$ being a $p$-polynomial of $z_i$. 

If we go back to the generators $x_1,\ldots,x_n$ of $L$ and the relators $u_1,\ldots,u_m$ of $G=L/J$ and apply the same transformations as we did to their images $y_1,\ldots,y_n$ and $v_1,\ldots, v_m$ then, according to what was mentioned about the automorphisms of the \frs before the proof of Lemma \ref{l001}, we obtain the desired generators and relations for $G$.
\end{proof}

Now we are ready to complete the proof of Theorem \ref{tbp}.

\begin{proof}
It follows from Lemma \ref{l001} that we can choose a presentation for $G=L/I$ in such a way that one of the generators, say $t$, is not involved  in the $p$-power portions of defining relations. We denote this generator by $t$ and other generators by $a_1,\ldots,a_{n-1}$. 

Let $M_k$  be a restricted ideal of $L$ generated by $a_1,\ldots,a_{n-1}$ and $\tp{k}$, where $k=1,2,\ldots$. Then it follows from \cite[2.7.4]{B} that $M_k$ is a \fr with free generators $a_{li}=(\mathrm{ad}\,t)^l(a_i)$,  $0\leq l <p^k$, $0\leq i\le n-1$ and $\tp{k}$. 

Now we consider an ideal $M$ of $L$ generated by $a_1,\ldots,a_{n-1}$ then this ideal is a \rl whose free generators are $a_{li}=(\mathrm{ad}\,t)^l(a_i)$,  $l=0,1,\ldots$, $0\leq i\le n-1$. This follows because $M=\cap_{k=1}^{\infty}M_k$ and each finite subset of the generating set of $M$ mentioned above is a subset of the free generating set of $M_k$, for an appropriate $k$. 

By Proposition \ref{p001}, $I\subset M$. The image of $M$ in $G$ is defined by the relations $w_{lj}=0$ where each $w_{lj}$ is a $p$-polynomial $(\mathrm{ad}\,t)^l(w_j)$, rewritten in terms of the generators $a_{lj}$, where $l$ is as above and $1\leq j \leq m$. Applying the  Leibniz rule for the derivations we can easily rewrite each $w_{j}$ in terms of a finite subset of the latter set of generators $a_{li}$. Thus me may assume that there is a number $q$ such that only $a_{li}$ with $l<p^{\,q}$ are involved in the expression of $w_j$ as the elements of $M$. 

Now let us choose $k$ so that $p^{\,q} < p^k$ and an element $s$, which is a $p$-polynomial in $t$ with leading term $\tp{k}$. Let $M(s)$ be an ideal of $L$ generated by $a_1,\ldots,a_{n-1}$ and $s$. Again by \cite[Section 2.7]{B}, $M(s)$ is freely generated by $a_{li}=(\mathrm{ad}\,t)^l(a_i)$,  $0\leq l <p^k$, $0\leq i\le n-1$ and $s$. The image $P$ of $M(s)$ in $G$ is defined by the relations $w_{lj}=(\mathrm{ad}\,t)^l(w_j)$, $l$ is as above, $1\leq j \leq m$, which have to be rewritten through the new free generating system. (This a known fact but it follows also from our Lemma \ref{l002} below). If we perform derivation in the relations $w_{lj}$ with the use of the Leibniz rule then we observe that the maximum value of index $l$ in the generators of $M$ that are involved in $w_{lj}$ is less than $p^{\,q}+p^k$. The generators $a_{li}$ of $M$ with $p^k\leq l< p^k+p^{\,q}$ are no longer on the above list of the free generators of $M(s)$. For these we have $a_{p^k+j,i}=(\mathrm{ad}\,t)^{p^k}(a_{ji})=[t^{[p^k]},a_{ji}]$, where $0\leq j < p^{\,q}$, and so if $s=\tp{k}+\sum_{m=0}^{k-1}\alpha_m\tp{m}$ for some $\alpha_m\in F$, then each $a_{p^k+j,i}$ should be replaced by $[s,a_{ji}]-\sum_{m=0}^{k-1}\alpha_m a_{j+p^m,\,i}$. Since $q,m<k$, each generator $a_{j+p^m,\,i}$ is on the list of the free generators of $M(s)$.

Now let us impose additional relations on $P$ by setting $a_{li}=0$ for all $i=1,\ldots,n-1$, $0\leq l<p^{\,q}$. Let $Q$ be the quotient algebra of $P$ obtained in this way. Then removing such superfluous generators we will be left with $(n-1)(p^k-p^{\,q})+1$ generators $s$ and $a_{li}$, $i$ as always, and $p^{\,q}\leq l<p^k$, and still $mp^k$ relations obtained by replacing some generators by 0. Since $s$ was involved in the relations of $P$ only inside the commutators $[s,a_{li}]$, where $0\leq l < p^{\,q}$, none of the newly obtained relations of $Q$ involves $s$. 

As a result, $Q$ is the free product of the subalgebra generated by $s$ and the subalgebra $K$ generated by $(n-1)(p^k-p^{\,q})$ generators $a_{li}$, $i$ as always, $p^{\,q}\leq l<p^k$, with $mp^k$ relations among them. The difference between the number of generators and relations for $K$ now takes the form $(n-m-1)p^k-(n-1)p^{\,q}$. It is now obvious that if we choose $k$ sufficiently large this latter number can be made positive and then by Lemma \ref{l001} $K$ can be mapped onto a free restricted algebra of rank 1. The free product $Q$ can then be mapped onto a free restricted algebra of rank 2. Therefore $P$ can be mapped onto a \fr of rank 2. It remains to notice that the codimension of $P$ in $G$ equals the codimension $k$ of $M(s)$ in $L$ and so is a finite number.
\end{proof}

\section{Lackenby - Olshanskii - Osin Theorem for Restricted Lie Algebras}\label{sLOOT}

Our aim in this section is the proof of Theorem \ref{tLOO}, which is an analogue of some  group-theoretical results in \cite{L} and \cite{OO} in the case of restricted Lie algebras.

Before we prove this theorem we need few lemmas.

\begin{lemma}\label{l002} Let $G$ be a Lie algebra, $N$ an ideal of $G$, $g$ an element of $N$. Let $C\subset C_G(g)$ where $C_G(g)$ is the centralizer of $g$ in $G$. Suppose that $T$ is any totally ordered subset of  $G$ whose union with $C+N$ spans $G$ as a vector space. Denote by $Z$ the set  
$$
\{(\ad t_1)\cdots (\ad t_k)(g)\vert\: t_1\leq\ldots\leq t_k\in T\}.
$$
Then the ideal of $G$ generated by $g$ coincides with the ideal of $N$ generated by $Z$. In the case where $G$ is a \rl over a field of characteristic $p >0$, and all the ideals are restricted, we can replace $Z$ by a subset $Z_p$ consisting of all monomial in which the degree of any $t_i$ is at most $p-1$.
\end{lemma}
\begin{proof} Since $N$ is an ideal of $G$, the above elements are in $N$. Now by the definition of the universal enveloping algebra $U(G)$, any ideal of $G$ is a left module for the adjoint representation of $U(G)$. Thus the ideal of $G$ generated by $g$ is a submodule of the left $U(G)$-module $G$ generated by $g$. Using PBW-theorem \cite[Chapter 5]{J}, if we choose a totally ordered basis of $G$, in which the elements $n_\alpha$ of $N$ precede some elements $t_\beta$ of $T$ and these precede some elements $c_\gamma$ of $C$, then any element of $U(G)$ is a linear combination of the ordered monomials of the form
$$
n_{\alpha_1}\cdots n_{\alpha_k}t_{\beta_1}\cdots t_{\beta_l}c_{\gamma_1}\cdots c_{\gamma_m}.
$$
The action of $U(G)$ on $G$ is the unique extension of the adjoint representation. If we apply the above monomial to $g$ and recall that $C$ is the centralizer of $g$ we will see that in $N$ the ideal in question is generated by the elements $(\ad t_{\beta_1})\cdots (\ad t_{\beta_l})(g)$ with $t_{\beta_1}\leq\ldots\leq t_{\beta_l}$, as claimed.

In the case where $G$ is a restricted Lie algebra over a field of characteristic $p>0$ the universal enveloping algebra should be replaced by the restricted enveloping algebra $u_p(G)$. Then the restricted ideals of $G$ are left $u_p(G)$-submodules closed under the $p$-operation. As mentioned in the Introduction, when we generate an ideal we can first apply the action of $u_p(G)$ and then take all possible $p$-powers. Thus the argument as just above applies also in this case. By PBW Theorem for restricted enveloping algebras \cite[Chapter 5]{J} any element $x_\alpha$, $t_\beta$, and $c_\gamma$ enters the monomials of the basis to the degree at most $p-1$, as claimed. Thus the proof is complete.
\end{proof}

Our next lemma is as follows.
\begin{lemma}\label{l003} For any finite collection of nonzero elements $g_1,\ldots,g_k$ of a \fr $L$ and any number $n\in\mathbb{N}$ there is $m\in \mathbb{N}$ with the following property. For every $q\geq m$ there is  a restricted ideal $N$ of finite codimension in $L$ such that for all $1\leq i\leq k$ we have $\Sp{g_i,g_i^{[p]},\ldots,g_i^{[p^{n-1}]}}\cap N=\{ 0\}$ but $g_i^{[p^{\,q}]}\in N$.
\end{lemma}
\begin{proof} It is sufficient to prove this lemma in the case where $L$ is finitely generated. If $x_1,\ldots,x_r$ is the basis of $L$ then there is a natural filtration on $L$ in which an element $w^{[p^k]}$ of the canonical basis has filtration $p^nd$ if $w$ is a commutator in $x_1,\ldots,x_r$ of degree $d$. An arbitrary $g\in L$ has filtration $f$ if $f$ is the least filtration of the basic elements in its expression through the basis. The set of elements of filtration at least $f$ is a restricted ideal of $L$ which we denote by $I_f$. Obviously each such ideal is of finite codimension in $L$. Now suppose $m_0$ is the maximum filtration of the elements $g_1,\ldots,g_k$. Choose $m=m_0p^n+1$. For each $q\geq m$, we set $N=I_m$. If $d_i$ is the filtration of $g_i$ then the filtration of $g_i^{[p^s]}$ equals $d_ip^s$. By our choice of $m$, for each $i$, $1\leq i\leq k$, the elements $g_i,g_i^{[p]},\ldots,g_i^{[p^{n-1}]}$ are linearly independent modulo $N$. But if $q\geq m$ then for any element $a\in L$ we always have $a^{[p^{\,q}]}\in I_m$.
\end{proof}

One more result we need for the proof of Theorem \ref{tLOO} is the following.

\begin{proposition}\label{l004} Let $L$ be a \fr of rank $r\geq 2$, $g_1$, \ldots, $g_k$ arbitrary elements of $L$. Let $J_q$ be a restricted ideal of $L$ generated by the elements $g_1^{[p^{\,q}]},\ldots ,g_k^{[p^{\,q}]}$, where $q$ is a natural number. Then $\overline{L}=L/J_q$ is large for all but finitely many $q\in\mathbb{N}$.
\end{proposition}

\begin{proof} Without loss of generality we may assume that all $g_1,\ld,g_k$ are nonzero. By Lemma \ref{l003}, there exists $m\in\nat$ such that for any $q\geq m$ there is a restricted ideal $N$ of finite codimension such that for all $i=1,\ld,k$ we have the elements $g_i,g_i^{[p\,]},\ldots,g_i^{[p^{k}]}$ are linearly independent modulo $N$, but $g_i^{[p^{\,q}]}\in N$. In particular, the codimension of $\sa{g_i}+N$ is bounded from above by $j-k-1$ where $\dim L/N=j$. Now we want to show that the image $\overline{N}$ of $N$ in $\overline{L}$ is a large algebra. Let $T_i$ be a minimal set of elements of $L$ such that the union of $T_i$ and $\sa{g_i}+N$ span $L$ as a vector space. Then according to Lemma \ref{l002}, $\overline{N}$ is isomorphic to the quotient algebra of $N$ by the restricted ideal generated by the elements of the set $Z$ where $Z=\bigcup_{i=1}^kZ_i$ and 
$$
Z_i=\{(\ad t_1)^{l_1}\cdots (\ad t_{s_i})^{l_{s_i}}(g_i^{[p^{\, q}]})\vert\: \{t_1,\ldots,t_{s_i}\}= T_i,\: 0\leq l_1,\ld,l_{s_i}<p\}.
$$
Now since $\dim L/N=j$, according to the analogue of Schreier's Formula \cite[Theorem 2.7.5]{B} for the number of generators of a subgroup of a free group, the number of generators of $N$ is  $p^j(r-1)+1\geq p^j+1$. Now the codimension $r_i$ of $\sa{g_i}+N$ is at most $j-k-1$ and so the number of elements in each $Z_i$ is at most $p^{j-k-1}$. In this case the total number of defining relations for $\overline{N}$ is at most $kp^{j-k-1}<p^{j-1}\le p^j-1$. The difference between the number of generators and relations will be at least 2 and so Theorem \ref{tbp} applies proving that $\overline{N}$ is a large \rlp Since $\overline{N}$ is of finite codimension in $\overline{L}$ this latter is a large \rlc as required.
\end{proof}

Now we can comment on Proposition \ref{p881}.

\begin{proof} If $k=1$ then we can adapt the previous proof to get a stronger result of Proposition \ref{r001}. If $N$ is an ideal of $L$ of codimension $j$, $g\notin N$ and $h\in\sa{g}\cap N$ then the same proof as just above shows that the number of generators of $N$ is still  $p^j(r-1)+1\geq p^j+1$ and $\da{\gp{n}}{L}$ is generated as an ideal of $N$ by $p^{j-1}$ elements. The difference is greater than $p^j+1-p^{j-1}\geq 2$. As before, we use Theorem \ref{tbp} to derive that $N/\da{h}{L}$ is large. This proves Proposition \ref{p881}, which is important in the proof of Theorem \ref{tOO}.\end{proof}

\smallskip

The following lemma shows that the situation in restricted Lie algebras can be very different from that in groups. In the case of groups, if we are given an element $g$ of a normal $H$ of index $d$ in a group $G$, the normal subgroup $\da{g}{G}$ of G generated by $g$ is a normal subgroup of $H$ generated by the conjugates $x_1gx_1^{-1},\ldots,x_dgx_d^{-1}\in H$ for some $x_1,\ldots, x_d\in G$. In the case of \rlsc we have instead the following.

\begin{lemma}\label{l991} Let $G$ be a finitely generated \rlc $H$ a restricted ideal of $G$ such that $\dim G/H=d$, $g\in H$, $n\in\nat$. Let $I=\da{g^{[p^n]}}{G}$ be the restricted ideal of $G$ generated by $g^{[p^n]}$, where $d\le n$, $J=\da{g^{[p^{n-d}]}}{H}$ the restricted ideal of $H$ generated by $g^{[p^{n-d}]}$. Then $I\subset J$.
\end{lemma}

\begin{proof}By Lemma \ref{l002}, $I$ as an ideal of $H$ is generated by the elements of the set defined as follows
\begin{equation}\label{e001}
Z=\{(\ad t_1)^{l_1}\cdots (\ad t_{d})^{l_{d}}(g^{[p^n]})\vert\: \{t_1,\ldots,t_{d}\}= T,\: 0\leq l_1,\ld,l_{d}<p\},
\end{equation}
where $T$ spans $G$ with $H$. Applying induction by $d$ with obvious basis for $d=0$, it is enough to show that $(\ad t)^i(\gp{n})\in\da{\gp{n-1}}{H}$, for $i=0,1,\ldots,p-1$. If $i>0$ then using the main identity of \rls (\ref{em}) and the Leibniz rule, one can write a  commutator formula as follows. All commutators are left-normed, that is $[u,v,w]=[u,[v,w]]$, for any $u,v,w\in G$.
\begin{eqnarray*}&&(\ad t)^i(\gp{n})=[\underbrace{t,\ldots,t}_{i},\gp{n}]= -[\underbrace{t,\ldots,t}_{i-1},[\gp{n},t]]\\&=&-[\underbrace{t,\ldots,t}_{i-1},[\underbrace{\gp{n-1},\ldots,\gp{n-1}}_{p},t]]\\&=&[\underbrace{t,\ldots,t}_{i-1},[\underbrace{\gp{n-1},\ldots,\gp{n-1}}_{p-1},[t,\gp{n-1}]]]
\end{eqnarray*}
Thus, applying the Leibniz rule, we can write 
\begin{eqnarray*}&&(\ad t)^i(\gp{n})\\&=&\sum_{k_1,\ldots k_p}[(\ad t)^{k_1}(\gp{n-1}),[\ldots,[(\ad t)^{k_{p-1}}(\gp{n-1}),(\ad t)^{k_p}(\gp{n-1})]]].
\end{eqnarray*} 
It is required, in the latter sum, that $k_1+\cdots+k_p=i$ and $k_p>0$. Now for each $j=1,\ldots,p$ one has
$(\ad t)^{k_j}(\gp{n-1})\in H$. Also, because $i\leq p$, some $k_j=0$. Thus, indeed, if $i>0$, the expression in question is in the ideal generated by $\gp{n-1}$. If $i=0$ then $(\ad t)^0(\gp{n})=(\gp{n-1})^{[p]}$ is in the restricted ideal generated by $\gp{n-1}$.
\end{proof}

Now we can complete the proof of Theorem \ref{tLOO}.

\begin{proof} 

As just proved, $I_n$ is contained in the ideal $J_{m}$ generated in $H$ as a restricted ideal by the elements $g_1^{[p^{m}]},\ldots,g_k^{[p^{m}]}$, where $m=n-d$ and $d=\dim G/H$. If we prove that the homomorphic image $H/J_{m}$ of $H/I_n$ is large, then also $H/I_n$ is large. Since $H/I_n$ is an ideal of finite codimension in $G/I_n$, we will be able to conclude that $G/I_n$ is a large \rlp 

Now let $\ve$ be a homomorphism of $H$ onto a nonabelian \fr $L$, assumed in the statement of our theorem. Set $K_m=\ve(J_m)$. Then $H/J_m$ admits a surjective homomorphism onto $L/K_m$.  Now $K_m$ is generated in $L$, as a restricted ideal, by the elements $\ve(g_1)^{[p^{m}]},\ldots,\ve(g_k)^{[p^{m}]}$. By Lemma \ref{l004}, for any set of elements $\ve(g_1),\ldots,\ve(g_k)\in L$, there is a number $M$ such that if $m > M$ and $K_m$ is the restricted ideal generated by $\ve(g_1)^{[p^{\,m}]},\ldots,\ve(g_k)^{[p^{\,m}]}$ then $L/K_m$ is large. Hence $H/J_m$ is large by Proposition \ref{ppla1}, Claim (i), as desired. Thus the proof is complete.
\end{proof}

\section{Constructing Nil Restricted Lie Algebras}\label{sCPLA}

In what follows we will use the \textit{derived $p$-series} $\{\de{i}{G}\vert\,i=0,1,\ld\}$ of derivation stable ideals of a \rl $G$ defined as follows. We set $\de{0}{G}=G$ and
$$
\de{i}{G}=[\de{i-1}{G},\de{i-1}{G}]+(\de{i-1}{G})^{[p\,]}\mbox{ for }i\ge 1.
$$
Obviously, $G/\de{1}{G}$ is finite-dimensional as a finitely generated abelian \rl with all elements of nil-index $p$. By \cite[2.7.5]{B} then $\de{1}{G}$ is a finitely generated \rlp Continuing in the same way, we obtain that each algebra $\de{i-1}{G}/\de{i}{G}$ is finite-dimensional and each $\de{i}{G}$ is finitely generated, for $i=1,2,\ldots$. Thus each algebra $G/\de{i}{G}$ is finite-dimensional and applying (\ref{em}) and Engel's Theorem \cite[1.7.3]{B} we easily derive that each $G/\de{i}{G}$ is nilpotent as a Lie algebra.

 The following Lemma is immediate using an argument similar to the one used in Lemma \ref{l003}.

\begin{lemma}\label{l005} For any finite-dimensional subspace $V$ of a \fr $L$ there exists $d\in\nat$ such that $\de{d}{L}\cap V=\{ 0\}$.
\end{lemma}

The next proposition is a version of Theorem \ref{tLOO}.

\begin{proposition}\label{p002} Let $G$ be a finitely generated \rlp Suppose that $P$ is an ideal of finite codimension in $G$ that can be mapped homomorphically onto a nonabelian \frp Then for any element $g\in P$ there is $m\in\nat$ such that if $\gp{n}\in \de{m}{G}$ then $\de{m}{P}/\da{\gp{n}}{G}$ can be mapped homomorphically onto a nonabelian \frp
\end{proposition}

\begin{proof} Let $a=\dim L/P$. We start with proving a weaker statement. 

\begin{claim}\label{c} Given an arbitrary $g\in P$ there exists $m$ such that if $\gp{n+a}\in \de{m}{G}$ then $\de{m}{G}/\da{\gp{n}}{P}$ can be mapped homomorphically onto a nonabelian \frp
\end{claim}

Let $\ve$ be a surjective homomorphism $\ve:P\rightarrow L$, $L$ a \fr and $h=\ve(g)$. If $h=0$ then we may set $m=1$. Suppose $g^{[p^n]}\in \de{1}{P}$ and $I=\da{\gp{n+a}}{G}$. Then by Lemma \ref{l991}, $I$ is contained in  $Q=\da{\gp{n}}{G}$. Thus $\de{1}{P}/I$ is naturally mapped onto $\de{1}{P}/Q$. Since under $\ve$ the element $g$ is mapped into 0, $\ve$ induces a homomorphism $\bar{\ve}$ of $P/Q$ onto $L$. Applying Proposition \ref{ppla2}, Part (i), we find that $\de{1}{P}/Q$ maps homomorphically onto a nonabelian \frp  So in the case where $h=0$, the proof is complete.

Thus we may assume that $h\neq 0$. By Lemma \ref{l005}, there is $d\in\nat$ such that $h\notin\de{d}{L}$. We set $M=\de{d}{L}$. Then $\hp{d}\in M$. If we denote by $J$ the ideal of $L$ generated by $\hp{d}$ then by Proposition \ref{p881}, $\overline{M}=M/J$ is a large Lie algebra with a presentation in which the number of generators exceeds the number of relations at least by 2. In this case, according to Proposition \ref{r001} after Theorem \ref{tbp}, there is a restricted ideal $\bN$ of $\bM$ such that $\bM/\bN$ is a nil cyclic \rl of nil-index $c$, for a natural number $c$,  and there is a surjective homomorphism $\eta: \bN\rightarrow L_1$ where $L_1$ is a nonabelian \frp\ It follows from the definition of the derived $p$-series that  $\de{c}{\bM}\subset\bN$. Since $\de{c}{\bM}$ is an ideal of finite codimension in $\bN$, we may apply Proposition \ref{ppla2}, Part (i), to derive that $\de{c}{\bM}$ can be mapped homomorphically on a \frp

Let us set $m=d+c$. Assume $\gp{n}\in\de{m}{P}$. Using Lemma \ref{l991}, we conclude that $I=\da{\gp{n+a}}{G}$ as a restricted ideal of $P$ is contained in the ideal $Q=\da{\gp{n}}{P}$. Thus $\de{m}{P}/I$ maps homomorphically onto $\de{m}{P}/Q$. If we apply a homomorphism induced by $\ve: P\rightarrow L$ then $\de{m}{P}/J$ will be mapped onto $\de{d+c}{L}/R=\de{c}{M}/R$ where $R=\da{\hp{n}}{L}$. Since $R\subset J$, there is a natural homomorphism of $\de{c}{M}/R$ onto $\de{c}{\bM}$ which is mapped onto a \frc as it was shown above. Thus we have established Claim \ref{c}.

Now we can complete the proof of Proposition \ref{p002}. Indeed, once $m$, as in Claim \ref{c}, has been found, choose $n_0$ minimal such that $\gp{n_0}\in \de{m}{P}$. Then choose a \textit{maximal} number $m'\ge m$ such that $\gp{n_0+a}\in \de{m'}{P}$. Note that $\gp{n_0+a-1}\notin \de{m'}{P}$. This $m'$ is the number $m$ sought for $g$ in Proposition \ref{p002}. For, if $\gp{n'}\in \de{m'}{P}$ then $\gp{n'-a}\in \de{m}{P}$ and so $\de{m}{P}/\da{\gp{n'}}{G}$ can be mapped homomorphically on a nonabelian \frp Since $\de{m'}{P}/\da{\gp{n'}}{G}$ is a subalgebra of finite codimension in  $\de{m}{P}/\da{\gp{n'}}{G}$,  by Proposition \ref{ppla1} (ii), it can be mapped homomorphically on a \frp
\end{proof}

Now we can give a construction of  infinite - dimensional finitely generated nil \rls and their generalizations claimed in Theorem \ref{tOO} and Corollaries \ref{cOO0}, \ref{c001} and \ref{c002}. This construction is an analogue of a group-theoretic one due to Olshanskii - Osin \cite{OO}.

\begin{proof} Let $G$ be a finitely generated \rl and $P$ an ideal that can be mapped homomorphically onto a \frp
 Let $\{f_1, f_2,\ldots\}$ be the list of all elements of $P$. We set $G_0=G$, $P_0=P$ and suppose we already constructed \rls $G_i\supset P_i$ which are the homomorphic images of $G$ and $P$ (and of $G_{i-1}$, $P_{i-1}$ for $i>0$), with the same kernel, so that $\de{r_i}{P_i}$ maps homomorphically on a nonabelian \fr for some $r_i\in\nat$ and in which all the images of $f_1,\ldots,f_i$ are nilpotent. Let $g_{i+1}$ denote the image of $f_{i+1}$ in $G_i$. Then we choose $r_{i+1}>r_i$ and $n=n(i)\in\nat$, according to Proposition \ref{p002}, so that $g^{[p^n]}_{i+1}\in\de{r_{i+1}}{P_i}$ and also $\de{r_{i+1}}{P_i}/\da{g^{[p^n]}_{i+1}}{G_i}$ maps homomorphically onto a nonabelian \frp We set $G_{i+1}=G_i/\da{g^{[p^n]}_{i+1}}{G_i}$ and $P_{i+1}=P_i/\da{g^{[p^n]}_{i+1}}{G_i}$. By Proposition \ref{ppla2}, Part (i), $\de{r_{i+1}}{P_{i+1}}$ maps homomorphically on a nonabelian \frp Since  $\de{r_i}{P_i}$  maps homomorphically on a nonabelian \fr and $r_{i+1}>r_i$, it easily follows that $\de{r_{i+1}}{P_i}$ is a proper subalgebra of $\de{r_{i}}{P_i}$. Therefore,
\begin{equation}\label{e881}
\dim P_{i+1}/\de{r_{i+1}}{P_{i+1}}=\dim P_{i}/\de{r_{i+1}}{P_{i}}> \dim P_{i}/\de{r_{i}}{P_{i}}.
\end{equation}

Let $K_i$ denote the common kernel of the natural homomorphisms $G\rightarrow G_i$ and $P\rightarrow P_i$. Clearly, $\bP=P/\bigcup_{i=0}^{\infty}K_i$ is a nil restricted ideal in the \rl $\bG=G/\bigcup_{i=0}^{\infty}K_i$ Further we set $\wP=\bP/\bigcap_{r=0}^{\infty}\de{r_i}{\bP}$, $\wG=\bG/\bigcap_{r=0}^{\infty}\de{r_i}{\bP}$. Then $\wP$ is nil and residually finite-dimensional nilpotent, $\wG$ is residually finite dimensional, and $\wG/\wP\cong G/P$. To show that $\wG$ is infinite-dimensional, we observe that $\mathrm{Ker}\,(P\rightarrow P_i)\subset\de{r_i}{P_i}$, for every $i$. Hence $P/\de{r_{i}}{P_{i}}\cong P_{i}/\de{r_{i}}{P_{i}}$. Now by (\ref{e881}), $\dim P_{i}/\de{r_{i}}{P_{i}}\rightarrow\infty$ as $i\rightarrow\infty$. Therefore $\wP$ is  infinite-dimensional since it maps homomorphically onto $P_{i}/\de{r_{i}}{P_{i}}$ for every~$i$.
\end{proof}

Now we derive the Corollaries \ref{cOO0}, \ref{c001} and \ref{c002}. 

\smallskip

In the case of Corollary \ref{cOO0}, we know from Proposition \ref{r001} that any \rl $G$ presented as described has an ideal $P$ such that $G/P$ is cyclic nil. If we apply to $G$ and $P$ the construction of the previous theorem, we obtain a finitely generated $\wG$ with a nil-ideal $\wP$ such that $\wG/\wP\cong G/P$ is nil. Obviously, then $\wG$ is itself nil. Also, by Engel's Theorem, any finite-dimensional nil Lie algebra is nilpotent. Whence our claim about $\wG$ being residually finite-dimensional nilpotent. Now Corollary \ref{c001} is a direct consequence of Corollary \ref{cOO0}.

As for Corollary \ref{c002}, we start with a finitely generated infinite - dimensional nil \rl $G$, over an algebraic closure $\overline{F}$ of $F$. Then we can consider an ordinary Lie algebra $\wg$ generated by a finite generating set $X$ of $G$ over $F$. As mentioned earlier, any element in $G$ is a linear combination with coefficients in $\overline{F}$ of $p$-powers of the Lie monomials in $X$, that is, the $p$-powers of the elements in $\wg$. Were $\wg$ finite-dimensional then using another basic identity of \rlsc $(x+y)^{[p]}=x^{[p]}+y^{[p]}+w(x,y)$ for any $x$, $y$ in $L$ and $w(x,y)$ in the ordinary Lie subring generated by $x,y$, we would easily obtain $G$ finite - dimensional (over $\overline{F}$). Notice, that by (\ref{em}) each $\mathrm{ad}\, g$ is nilpotent. So $\wg$ is an example of an infinite-dimensional finitely generated Engel Lie algebra.

\end{document}